\documentclass[12pt]{article}
\usepackage{amsmath, amsthm, amssymb}
\title{Generic canonical form of pairs
of matrices with zeros\footnotetext{This is the authors' version of a work that was published in Linear Algebra Appl. 380 (2004) 241-251.}}
\author{Tat$'$yana N. Gaiduk%
\\Department of Physics and Mathematics\\
Chernigov Pedagogical University,
Chernigov, Ukraine
 \and
Vladimir V. Sergeichuk%
\thanks{Corresponding author.
Partially supported by NSF grant
DMS-0070503. {\it E-mail address:}
sergeich@imath.kiev.ua}
\\ Institute of
Mathematics\\ Tereshchenkivska 3, Kiev,
Ukraine}
\date{}

\begin{document}

\renewcommand{\le}{\leqslant}
\renewcommand{\ge}{\geqslant}
\newcommand{\rank}%
  {\mathop{\rm rank}\nolimits}
\newcommand{\size}%
  {\mathop{\rm size}\nolimits}
\newcommand{\cd}%
  {\text{\raisebox{3pt}{$\centerdot$}}}
\newcommand{\mi}%
  {^{\scriptscriptstyle -}}
\newcommand{\pl}%
  {^{\scriptscriptstyle +}}
\newcommand{\lin}%
  {\frac{\quad}{}}

\newtheorem{theorem}{Theorem}
\newtheorem{lemma}[theorem]{Lemma}

\maketitle
\begin{abstract}
We consider a family of pairs of
$m\times p$ and $m\times q$ matrices,
in which some entries are required to
be zero and the others are arbitrary,
with respect to transformations
$(A,B)\mapsto (SAR_1,SBR_2)$ with
nonsingular $S,\ R_1$, and $R_2$. We
prove that almost all of these pairs
reduce to the same pair $(A_0, B_0)$
from this family, except for pairs
whose arbitrary entries are zeros of a
certain polynomial. The polynomial and
the pair $(A_0, B_0)$ are constructed
by a combinatorial method based on
properties of a certain graph.

{\it AMS classification:} 15A21

{\it Keywords:} Structured matrices;
Parametric matrices; Canonical forms
\end{abstract}

\section{Introduction and main results}
\label{s1}

Let ${\cal A}:U_1\to V$ and ${\cal
B}:U_2\to V$ be linear mappings of
vector spaces over an arbitrary field
$\mathbb F$. Changing the bases of the
vector spaces, we may reduce the
matrices $A$ and $B$ of these mappings
by transformations
\begin{equation}\label{1.1}
 (A,B)\mapsto (SAR_1,SBR_2)
 \qquad\text{with nonsingular
 $S,\ R_1$, and $R_2$.}
\end{equation}
A canonical form of $(A,B)$ for these
transformations is
\begin{equation}\label{1.eee}
\left(\begin{bmatrix}
   I_r&0&0\\
   0&I_s&0\\
   0&0&0  \\
   0&0&0
\end{bmatrix},\:
\begin{bmatrix}
 0&I_r& 0\\
 0&0&0\\
I_t&0&0\\ 0&0&0
\end{bmatrix}
\right),
\end{equation}
where $I_r$ denotes the $r$-by-$r$
identity matrix and $r$, $s$, and $t$
are determined by the equalities
$r+s=\rank A$, $r+t=\rank B$, and
$r+s+t=\rank\, [A\,|\,B]$ (see Lemma
\ref{t1}).

We consider a family of pairs $(A,B)$,
in which $n$ entries $a_1,\dots,a_n$
are arbitrary and the others are
required to be zero. We prove that
there exists a nonzero polynomial
$f(x_1,\dots,x_n)$ such that all pairs
$(A,B)$ with $f(a_1,\dots,a_n)\ne 0$
reduce to the same pair $(A_{\text{\rm
gen}},B_{\text{\rm gen}})$ from this
family. The pair $(A_{\text{\rm
gen}},B_{\text{\rm gen}})$ has the form
\eqref{1.eee} up to permutations of
columns and simultaneous permutations
of rows in $A$ and $B$. Following
\cite{wou}, we call $(A_{\text{\rm
gen}},B_{\text{\rm gen}})$ a
\emph{generic canonical form} of the
family (this notion has no sense if
$\mathbb F$ is a finite field). We give
a combinatorial method of finding
$f(x_1,\dots,x_n)$ and $(A_{\text{\rm
gen}},B_{\text{\rm gen}})$.

\subsection{Generic canonical form
of matrices with zeros}

Since the rows of $A$ and $B$ in
\eqref{1.1} are transformed by the same
matrix $S$, we represent the pair
$(A,B)$ by the block matrix
$
M=[A\,|\,B],
$
which will be called a {\it bipartite
matrix}. A family of bipartite
matrices, in which some entries are
zero and the others are arbitrary, may
be given by a matrix
\begin{equation}\label{1.7}
M(x)=[A(x)\,|\, B(x)], \quad
x=(x_1,\dots,x_n),
\end{equation}
whose $n$ entries are unknowns
$x_1,\dots,x_n$ and the others are
zero. For instance,
\begin{equation}\label{1.5}
M(x)= \left[\begin{array}{cc|ccc}
  0&0           &  x_4&x_7&0\\
  x_1&0   &  x_5&0&0\\
  0&x_2   &  0&0&x_9\\
  0&x_3   &  x_6&x_8&0
\end{array}\right]
\end{equation}
gives the family $\{M(a)\,|\,
a\in{\mathbb F}^9\}$.

Considering \eqref{1.7} as a matrix
over the field
\begin{equation}\label{3.2}
 {\mathbb K}=
 \left.\left\{
 \frac{f(x_1,\dots,x_n)}
 {g(x_1,\dots,x_n)}\,\right|
 \,f,g\in
 {\mathbb F}[x_1,\dots,x_n]
\text{ and } g\ne 0
 \right\}
\end{equation}
of rational functions (its elements are
quotients of polynomials), we put
\begin{equation}\label{1.7v}
r_A=\rank_{\mathbb K}A({x}),\quad
r_B=\rank_{\mathbb K} B({x}),\quad
r_M=\rank_{\mathbb K} M({x}).
\end{equation}

The following theorem is proved in
Section \ref{s_proof}.

\begin{theorem}\label{t0.1}
Let $M(x)=[A({x})\,|\, B({x})]$ be a
matrix whose $n$ entries are unknowns
$x_1,\dots,x_n$ and the others are
zero. Then there exists a nonzero
polynomial
\begin{equation}\label{1.7a}
  f(x)
  =\sum c_i x_1^{m_{i1}}\cdots
  x_n^{m_{in}}
\end{equation}
such that all matrices of the family
\begin{equation*}\label{n1}
  {\cal M}_f=\{M(a)\,|\,a\in{\mathbb
  F}^n\text{ and }f(a)\ne 0\}
\end{equation*}
reduce by transformations
$[A\,|\,B]\mapsto [SAR_1\,|\,SBR_2]$
with nonsingular $S,\ R_1$, and $R_2$
to the same matrix
\begin{equation}\label{1.6r}
M_{\text{\rm gen}} =[A_{\text{\rm
gen}}\,|\,B_{\text{\rm gen}}] \in{\cal
M}_f.
\end{equation}
Up to a permutation of columns within
$A_{\text{\rm gen}}$ and $B_{\text{\rm
gen}}$ and a permutation of rows, the
matrix $M_{\text{\rm gen}}$ has the
form
\begin{equation}\label{1.2}
\left[\begin{array}{ccc|ccc}
   I_r&0&0 &  0&I_r& 0\\
   0&I_s&0 &  0&0&0\\
   0&0&0   &  I_t&0&0\\
   0&0&0   &  0&0&0
\end{array}\right],
\end{equation}
which is uniquely determined by $M(x)$
due to the equalities
\begin{equation}\label{1.2a}
r+s=r_A,\quad r+t=r_B,\quad r+s+t=r_M
\quad(\text{see \eqref{1.7v}}).
\end{equation}
\end{theorem}

We call $M_{\text{\rm gen}}$ a {\it
generic canonical form} of the family
$\{M(a)\,|\,a\in{\mathbb F}^n\}$
because $M(a)$ reduces to $M_{\text{\rm
gen}}$ for all $a\in {\mathbb F}^n$
except for those in the proper
algebraic variety
$
\{a\in {\mathbb F}^n\,|\, f(a)=0\}.
$

\subsection{A combinatorial
method} \label{ss12}

The polynomial $f(x)$ and the matrix
$M_{\text{\rm gen}}$ can be constructed
by a combinatorial method: we represent
the matrix $M(x)=[A({x})\,|\, B({x})]$
by a graph and study its subgraphs.
Similar methods were applied in
\cite{her,rob1,rob2,wou} to square
matrices up to similarity and to
pencils of matrices.

The graph is defined as follows. Its
vertices are
\[
1,\dots,m,\,1\!\mi,\dots,
p\mi,\,1\!\pl,\dots, q\pl,
\]
where $m\times p$ and $m\times q$ are
the sizes of $A({x})$ and $B({x})$. Its
edges
\begin{equation}\label{e3}
\alpha_1,\dots,\alpha_n
\end{equation}
are determined by the unknowns
$x_1,\dots,x_n$: if $x_l$ is the
$(i,j)$ entry of $A({x})$ then
$\alpha_l:i\lin j\mi$ (that is,
$\alpha_l$ links the vertices $i$ and
$j\mi$), and if $x_l$ is the $(i,j)$
entry of $B({x})$ then $\alpha_l:i\lin
j\pl$. The edges between
$\{1,\dots,m\}$ and $\{1\!\mi,\dots,
p\mi\}$ are called \emph{left edges},
and the edges between $\{1,\dots,m\}$
and $\{1\!\pl,\dots, q\pl\}$ are called
\emph{right edges}.

For example, the matrix \eqref{1.5} is
represented by the graph\\[7pt]
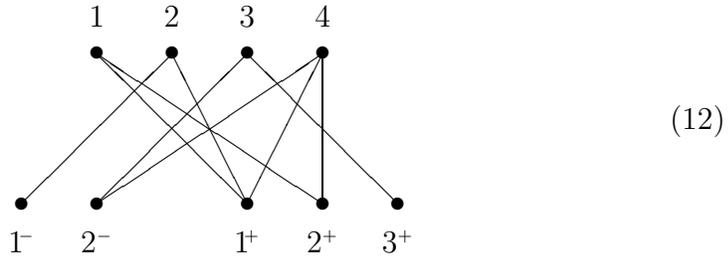
\begin{equation}\label{1.6}
\unitlength 1mm \linethickness{0.4pt}
\begin{picture}(50.00,15.00)
\put(0.00,-10.00){\makebox(0,0)[cc]{$\bullet$}}
\put(10.00,-10.00){\makebox(0,0)[cc]{$\bullet$}}
\put(30.00,-10.00){\makebox(0,0)[cc]{$\bullet$}}
\put(40.00,-10.00){\makebox(0,0)[cc]{$\bullet$}}
\put(50.00,-10.00){\makebox(0,0)[cc]{$\bullet$}}
\put(10.00,10.00){\makebox(0,0)[cc]{$\bullet$}}
\put(20.00,10.00){\makebox(0,0)[cc]{$\bullet$}}
\put(30.00,10.00){\makebox(0,0)[cc]{$\bullet$}}
\put(40.00,10.00){\makebox(0,0)[cc]{$\bullet$}}
\put(0.00,-15.00){\makebox(0,0)[cc]{$1\!\mi$}}
\put(10.00,-15.00){\makebox(0,0)[cc]{$2\mi$}}
\put(30.00,-15.00){\makebox(0,0)[cc]{$1\!\pl$}}
\put(40.00,-15.00){\makebox(0,0)[cc]{$2\pl$}}
\put(50.00,-15.00){\makebox(0,0)[cc]{$3\pl$}}
\put(10.00,15.00){\makebox(0,0)[cc]{$1$}}
\put(20.00,15.00){\makebox(0,0)[cc]{$2$}}
\put(30.00,15.00){\makebox(0,0)[cc]{$3$}}
\put(40.00,15.00){\makebox(0,0)[cc]{$4$}}
\put(0.00,-10.00){\line(1,1){20.00}}
\put(10.00,-10.00){\line(1,1){20.00}}
\put(10.00,-10.00){\line(3,2){30.00}}
\put(30.00,-10.00){\line(-1,1){20.00}}
\put(30.00,-10.00){\line(-1,2){10.00}}
\put(30.00,-10.00){\line(1,2){10.00}}
\put(40.00,-10.00){\line(-3,2){30.00}}
\put(40.00,-10.00){\line(0,1){20.00}}
\put(50.00,-10.00){\line(-1,1){20.00}}
\end{picture}
\\*[50pt]
\end{equation}
with the left edges
$\alpha_1,\,\alpha_2,\,\alpha_3$ and
the right edges
$\alpha_4,\,\alpha_5,\dots,\alpha_9$.

Each subset $\cal S$ in the set of
edges \eqref{e3} can be given by the
\emph{characteristic vector}
\begin{equation*}\label{4.03}
\varepsilon_{\cal
S}=(e_1,\dots,e_n),\qquad e_l=
  \begin{cases}
    1 & \text{if $\alpha_l\in\cal S$}, \\
    0 & \text{otherwise}.
  \end{cases}
\end{equation*}

By a {\it matchbox} we mean a set of
edges (={\it matches}) that have no
common vertices. The {\it size} of a
matchbox ${\cal S}$ is the number of
its matches; since each row and each
column of $M(\varepsilon_{\,\cal S})$
have at most one $1$ and the other
entries are zero,
\begin{equation}\label{1.6x}
\size {\cal S}= \rank
M(\varepsilon_{\,\cal S}).
\end{equation}
A matchbox is {\it left} ({\it right})
if all its matches are {left}
({right}). Such a matchbox is said to
be \emph{largest} if it has the maximal
size among all left (right) matchboxes.
For example, the subgraph\\[5pt]
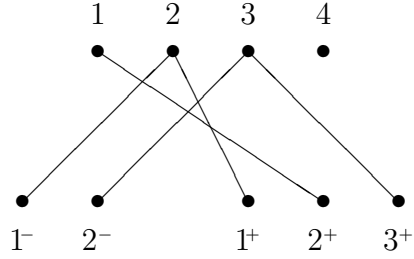
\begin{equation*}\label{1.6a}
\unitlength 1.00mm
\linethickness{0.4pt}
\begin{picture}(50.00,15.00)
\put(0.00,-10.00){\makebox(0,0)[cc]{$\bullet$}}
\put(10.00,-10.00){\makebox(0,0)[cc]{$\bullet$}}
\put(30.00,-10.00){\makebox(0,0)[cc]{$\bullet$}}
\put(40.00,-10.00){\makebox(0,0)[cc]{$\bullet$}}
\put(50.00,-10.00){\makebox(0,0)[cc]{$\bullet$}}
\put(10.00,10.00){\makebox(0,0)[cc]{$\bullet$}}
\put(20.00,10.00){\makebox(0,0)[cc]{$\bullet$}}
\put(30.00,10.00){\makebox(0,0)[cc]{$\bullet$}}
\put(40.00,10.00){\makebox(0,0)[cc]{$\bullet$}}
\put(0.00,-15.00){\makebox(0,0)[cc]{$1\!\mi$}}
\put(10.00,-15.00){\makebox(0,0)[cc]{$2\mi$}}
\put(30.00,-15.00){\makebox(0,0)[cc]{$1\!\pl$}}
\put(40.00,-15.00){\makebox(0,0)[cc]{$2\pl$}}
\put(50.00,-15.00){\makebox(0,0)[cc]{$3\pl$}}
\put(10.00,15.00){\makebox(0,0)[cc]{$1$}}
\put(20.00,15.00){\makebox(0,0)[cc]{$2$}}
\put(30.00,15.00){\makebox(0,0)[cc]{$3$}}
\put(40.00,15.00){\makebox(0,0)[cc]{$4$}}
\put(0.00,-10.00){\line(1,1){20.00}}
\put(10.00,-10.00){\line(1,1){20.00}}
\put(30.00,-10.00){\line(-1,2){10.00}}
\put(50.00,-10.00){\line(-1,1){20.00}}
\put(10.00,10.00){\line(3,-2){30.00}}
\end{picture}
\\*[50pt]
\end{equation*}
of \eqref{1.6} is formed by the largest
left and right matchboxes
\begin{equation}\label{1.6v}
{\cal A}=\{2\lin 1\!\mi,\ 3\lin
2\mi\}\quad\text{and}\quad {\cal
B}=\{1\lin 2\pl,\ 2\lin 1\!\pl,\ 3\lin
3\pl\}.
\end{equation}

For a left matchbox ${\cal A}$ and a
right matchbox ${\cal B}$, we denote by
$
{\cal A}\Cup{\cal B}
$
the matchbox obtained from ${\cal
A}\cup{\cal B}$ by removing all matches
of ${\cal B}$ that have common vertices
with matches of ${\cal A}$. For
example,
\begin{equation}\label{1.7w}
{\cal A}\Cup{\cal B}= \{2\lin 1\!\mi,\
3\lin 2\mi,\ 1\lin 2\pl\}
\end{equation}
for the matchboxes \eqref{1.6v}.

For every matchbox
\begin{equation*}\label{1.6m}
{\cal S}=\{i_1\lin {j_1\mi},\dots,\
i_{\alpha}\lin j_{\alpha}\mi,\
i_{{\alpha}+1}\lin k_{1}\pl,\dots,\
i_{\alpha+\beta}\lin k_{\beta}\pl\},
\end{equation*}
we denote by $\mu_{\cal S}({x})$ the
minor of order $\alpha+\beta$ in
$M({x})=[A({x})\,|\, B({x})]$ whose
matrix belongs to the rows numbered
$i_1,\dots,i_{\alpha+\beta}$, to the
columns of $A({x})$ numbered
${j_1},\dots,{j_{\alpha}}$, and to the
columns of $B({x})$ numbered
${k_1},\dots,k_{\beta}$. For example,
the matchbox \eqref{1.7w} determines
the minor
\[
\mu_{{\cal A}\Cup{\cal B}}(x)=
\begin{vmatrix}
 0&0           & x_7\\
 x_1&0   & 0\\
 0&x_2   & 0
\end{vmatrix}=x_1
x_2x_7 \quad\text{in \eqref{1.5}.}
\]

The next theorem will be proved in
Section \ref{s_proof}.

\begin{theorem}\label{t1.2}
The generic canonical form
$M_{\text{\rm gen}}$ and the polynomial
$f(x)$ from Theorem \ref{t0.1} may be
constructed as follows. We represent
$M(x)$ by the graph. Among pairs
consisting of a largest left matchbox
and a largest right matchbox, we choose
a pair $({\cal A},{\cal B})$ with the
minimal number $v({\cal A},{\cal B})$
of common vertices, and take
\begin{equation}\label{1.6dd}
M_{\text{\rm gen}}=
M(\varepsilon_{{\cal A}\cup{\cal B}}),
\qquad f(x)=f_{{\cal A}{\cal B}}(x),
\end{equation}
where $f_{{\cal A}{\cal B}}(x)$ is the
lowest common multiple of $\mu_{\cal
A}(x)$, $\mu_{\cal B}(x)$, and
$\mu_{{\cal A}\Cup{\cal B}}(x)$:
\begin{equation}\label{e6}
f_{{\cal A}{\cal B}}(x)= {\rm
LCM}\{\mu_{\cal A}(x),\ \mu_{\cal
B}(x),\ \mu_{{\cal A}\Cup{\cal
B}}(x)\}.
\end{equation}
Up to permutations of columns within
$A_{\text{\rm gen}}$ and $B_{\text{\rm
gen}}$ and a permutation of rows, the
matrix $M(\varepsilon_{{\cal
A}\cup{\cal B}})$ has the form
\eqref{1.2} with
\begin{equation}\label{1.6p}
r=v({\cal A},{\cal B}),\quad s=\size
{\cal A}-r, \quad\text{and}\quad
t=\size {\cal B}-r.
\end{equation}
\end{theorem}

\subsection{An example}

Let us apply Theorems \ref{t0.1} and
\ref{t1.2} to the family given by the
matrix \eqref{1.5} with the graph
\eqref{1.6}. The matchboxes
\eqref{1.6v} do not satisfy the
conditions of Theorem \ref{t1.2}
because they have two common vertices
`2' and `3'. This number is not minimal
since the largest matchboxes
\begin{equation}\label{e2}
{\cal A} = \{2\lin 1\!\mi,\ 3\lin
2\mi\},\qquad
 {\cal B} =
\{1\lin 1\!\pl,\ 3\lin 3\pl,\ 4\lin
2\pl\}
\end{equation}
forming the graph\\
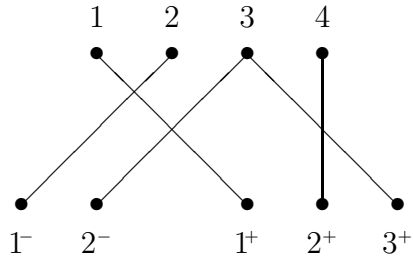
\begin{equation*}
\unitlength 1.00mm
\linethickness{0.4pt}
\begin{picture}(50.00,15.00)
\put(0.00,-10.00){\makebox(0,0)[cc]{$\bullet$}}
\put(10.00,-10.00){\makebox(0,0)[cc]{$\bullet$}}
\put(30.00,-10.00){\makebox(0,0)[cc]{$\bullet$}}
\put(40.00,-10.00){\makebox(0,0)[cc]{$\bullet$}}
\put(50.00,-10.00){\makebox(0,0)[cc]{$\bullet$}}
\put(10.00,10.00){\makebox(0,0)[cc]{$\bullet$}}
\put(20.00,10.00){\makebox(0,0)[cc]{$\bullet$}}
\put(30.00,10.00){\makebox(0,0)[cc]{$\bullet$}}
\put(40.00,10.00){\makebox(0,0)[cc]{$\bullet$}}
\put(0.00,-15.00){\makebox(0,0)[cc]{$1\!\mi$}}
\put(10.00,-15.00){\makebox(0,0)[cc]{$2\mi$}}
\put(30.00,-15.00){\makebox(0,0)[cc]{$1\!\pl$}}
\put(40.00,-15.00){\makebox(0,0)[cc]{$2\pl$}}
\put(50.00,-15.00){\makebox(0,0)[cc]{$3\pl$}}
\put(10.00,15.00){\makebox(0,0)[cc]{$1$}}
\put(20.00,15.00){\makebox(0,0)[cc]{$2$}}
\put(30.00,15.00){\makebox(0,0)[cc]{$3$}}
\put(40.00,15.00){\makebox(0,0)[cc]{$4$}}
\put(0.00,-10.00){\line(1,1){20.00}}
\put(10.00,-10.00){\line(1,1){20.00}}
\put(50.00,-10.00){\line(-1,1){20.00}}
\put(30.00,-10.00){\line(-1,1){20.00}}
\put(40.00,-10.00){\line(0,1){20.00}}
\end{picture}\\*[50pt]
\end{equation*}
have a single common vertex `3'. The
matchboxes \eqref{e2} satisfy the
conditions of Theorem \ref{t1.2} since
there is no pair of largest matchboxes
without common vertices.

The conditions of Theorem \ref{t1.2}
also hold for the largest matchboxes
\[
{\cal A}'= \{2\lin 1\!\mi,\ 4\lin
2\mi\},\qquad
 {\cal B}'=
\{1\lin 2\pl,\ 2\lin 1\!\pl,\ 3\lin
3\pl\}
\]
forming the graph\\
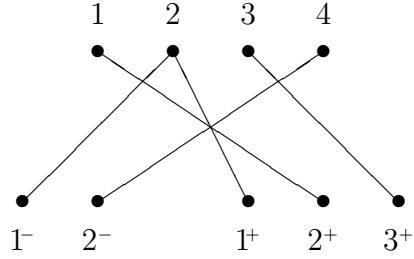
\begin{equation*}\label{1.8}
\unitlength 1.00mm
\linethickness{0.4pt}
\begin{picture}(50.00,15.00)
\put(0.00,-10.00){\makebox(0,0)[cc]{$\bullet$}}
\put(10.00,-10.00){\makebox(0,0)[cc]{$\bullet$}}
\put(30.00,-10.00){\makebox(0,0)[cc]{$\bullet$}}
\put(40.00,-10.00){\makebox(0,0)[cc]{$\bullet$}}
\put(50.00,-10.00){\makebox(0,0)[cc]{$\bullet$}}
\put(10.00,10.00){\makebox(0,0)[cc]{$\bullet$}}
\put(20.00,10.00){\makebox(0,0)[cc]{$\bullet$}}
\put(30.00,10.00){\makebox(0,0)[cc]{$\bullet$}}
\put(40.00,10.00){\makebox(0,0)[cc]{$\bullet$}}
\put(0.00,-15.00){\makebox(0,0)[cc]{$1\!\mi$}}
\put(10.00,-15.00){\makebox(0,0)[cc]{$2\mi$}}
\put(30.00,-15.00){\makebox(0,0)[cc]{$1\!\pl$}}
\put(40.00,-15.00){\makebox(0,0)[cc]{$2\pl$}}
\put(50.00,-15.00){\makebox(0,0)[cc]{$3\pl$}}
\put(10.00,15.00){\makebox(0,0)[cc]{$1$}}
\put(20.00,15.00){\makebox(0,0)[cc]{$2$}}
\put(30.00,15.00){\makebox(0,0)[cc]{$3$}}
\put(40.00,15.00){\makebox(0,0)[cc]{$4$}}
\put(0.00,-10.00){\line(1,1){20.00}}
\put(50.00,-10.00){\line(-1,1){20.00}}
\put(10.00,-10.00){\line(3,2){30.00}}
\put(30.00,-10.00){\line(-1,2){10.0}}
\put(40.00,-10.00){\line(-3,2){30.00}}
\end{picture}\\*[50pt]
\end{equation*}
since they have a single common vertex
too.

For these pairs of matchboxes, we have
\begin{gather*}
{\cal A}\Cup{\cal B} = \{2\lin 1\!\mi,\
3\lin 2\mi,\ 1\lin 1\!\pl,\ 4\lin
2\pl\},
  \\
f_{{\cal A}{\cal B}}(x)= {\rm
LCM}\{x_1x_2,\ x_9(x_6x_7-x_4x_8),\
x_1x_2(x_4x_8-x_6x_7)\}
\end{gather*}
and
\begin{gather*}
{\cal A}'\Cup{\cal B}'= \{2\lin
1\!\mi,\ 4\lin 2\mi,\ 1\lin 2\pl,\
3\lin 3\pl\},
  \\
f_{{\cal A}'{\cal B}'}(x)= {\rm
LCM}\{x_1x_3,\ -x_5x_7x_9,\
-x_1x_3x_7x_9\}.
\end{gather*}

By Theorems \ref{t0.1} and \ref{t1.2},
\[
\left[\begin{array}{cc|ccc}
  0&0           &  a_4&a_7&0\\
  a_1&0   &  a_5&0&0\\
  0&a_2   &  0&0&a_9\\
  0&a_3   &  a_6&a_8&0\\
\end{array}\right]\quad
\text{with }a_1,\dots,a_9\in{\mathbb F}
\]
(see \eqref{1.5}) reduces to the matrix
\[
M(\varepsilon_{{\cal A}\cup{\cal B}})
=\left[\begin{array}{cc|ccc}
  0&0           &  1&0&0\\
  1&0   &  0&0&0\\
  0&1   &  0&0&1\\
  0&0   &  0&1&0\\
\end{array}\right]
\ \text{if}\ f_{{\cal A}{\cal B}}(a)=
a_1a_2a_9(a_4a_8-a_6a_7)\ne 0
\]
and to the matrix
\[
M(\varepsilon_{{\cal A}'\cup{\cal
B}'})= \left[\begin{array}{cc|ccc}
  0&0   &  0&1&0\\
  1&0   &  1&0&0\\
  0&0   &  0&0&1\\
  0&1   &  0&0&0\\
\end{array}\right]
\ \text{if}\ f_{{\cal A}'{\cal B}'}(a)=
a_1a_3a_5a_7a_9\ne 0.
\]
Up to permutations of columns within
vertical strips and permutations of
rows, these matrices have the form
\[
\left[\begin{array}{cc|ccc}
  1&0   &  0&0&1\\
  0&1   &  0&0&0\\
  0&0   &  1&0&0\\
  0&0   &  0&1&0\\
\end{array}\right]\quad \text{(see
\eqref{1.2}).}
\]

\section{Proof of Theorems
\ref{t0.1} and \ref{t1.2}}
\label{s_proof}

\subsection{Bipartite matrices}
\label{s2}

The canonical form of a pair for
transformations \eqref{1.1} is well
known, see \cite[Sect. 1.2]{gab-roi}.
We recall it since we will use it in
the proof of Theorems \ref{t0.1} and
\ref{t1.2}.

Clearly, $(A,B)$ reduces to $(A',B')$
by transformations \eqref{1.1} if and
only if $[A\,|\,B]$ reduces to
$[A'\,|\,B']$ by a sequence of
\begin{itemize}
  \item[(i)]
elementary row-transformations in
$[A\,|\,B]$,
  \item[(ii)]
elementary column-transformations in
$A$, and
  \item[(iii)]
elementary column-transformations in
$B$.
\end{itemize}

\begin{lemma}\label{t1}
Every bipartite matrix $M=[A\,|\,B]$
over a field\/ $\mathbb F$ reduces by
transformations {\rm(i)--(iii)} to the
form
\begin{equation}\label{1.2'}
\left[\begin{array}{ccc|ccc}
   I_r&0&0 &  0&I_r& 0\\
   0&I_s&0 &  0&0&0\\
   0&0&0   &  I_t&0&0\\
   0&0&0   &  0&0&0
\end{array}\right]
\end{equation}
determined by the equalities
\begin{equation}\label{1.3}
r+s=\rank A,\quad r+t=\rank B,\quad
r+s+t=\rank M.
\end{equation}
\end{lemma}

\begin{proof}
By transformations (i) and (ii), we
reduce $M$ to the form
\begin{equation*}\label{2.2}
\left[\begin{array}{cc|c}
  I_h&0&B_1\\
  0&0&B_2
\end{array}\right],
\end{equation*}
and then by elementary
row-transformations within the second
horizontal strip and by transformations
(iii) to the form
\begin{equation*}\label{2.3}
\left[\begin{array}{cc|cc}
  I_h&0&B_3&B_4\\
  0&0&I_t&0\\
  0&0&0&0
\end{array}\right].
\end{equation*}
Adding linear combinations of rows of
$I_t$ to rows of $B_3$ by
transformations (i), we ``kill'' all
non-zero entries of $B_3$ :
\begin{equation*}\label{2.4}
\left[\begin{array}{cc|cc}
  I_h&0&0&B_4\\
  0&0&I_t&0\\
  0&0&0&0
\end{array}\right].
\end{equation*}
At last, we reduce $B_4$ to $I_r\oplus
0$ by elementary transformations. The
row-transformations with $B_4$ have
``spoiled'' the block $I_h$, but we
restore it by the inverse
column-transformations (ii) and obtain
the matrix \eqref{1.2'} with $r+s=h$.

Since the transformations (i)--(iii)
with $M=[A\,|\,B]$ preserve the ranks
of $M$, $A$, and $B$, we have the
equalities \eqref{1.3}. This implies
the uniqueness of \eqref{1.2'} since
$s=\rank M-\rank B$, $t=\rank M-\rank
A$, and $r=\rank A+\rank B-\rank M$.
\end{proof}

\subsection{Reduction of bipartite matrices
by permutations of rows and columns}
\label{s2a}

In this section we consider a bipartite
matrix $M=[A\,|\,B]$ with respect to
permutations of rows and columns.

\begin{lemma}\label{l3.1}
Every bipartite matrix $[A\,|\,B]$ with
linearly independent columns reduces by
a permutation of rows to the form
\begin{equation}\label{3.01}
\left[\begin{array}{c|c}
  A'     & \cd\\
  \cd&B'\\
  \cd&\cd
\end{array}\right],
\end{equation}
where $A'$ and $B'$ are nonsingular
square blocks and the points denote
unspecified blocks.
\end{lemma}

\begin{proof}
By permutations of rows we reduce
$[A\,|\,B]$ to the form
\begin{equation*}\label{3.02}
\left[\begin{array}{c|c}
  A_1&B_1\\
  \cd&\cd
\end{array}\right]
\end{equation*}
with a nonsingular square matrix
$[A_1\,|\,B_1]$. Laplace's theorem (see
\cite[Theorem 2.4.1]{pra}) states that
the determinant of $[A_1\,|\,B_1]$ is
equal to the sum of products of the
minors whose matrices belong to the
rows of $A_1$ by their cofactors
(belonging to $B_1$). One of these
summands is nonzero since
$[A_1\,|\,B_1]$ is nonsingular. We
collect the rows of the minor from this
summand at the top and obtain the
matrix \eqref{3.01}.
\end{proof}

\begin{lemma}\label{l3.2}
Every bipartite matrix $[A\,|\,B]$
reduces by permutations of rows and
permutations of columns in $A$ and $B$
to the form
\begin{equation}\label{3.03}
\left[\begin{array}{ccc|ccc}
  X_r&\cd&\cd     & \cd&Y_r&\cd\\
  \cd&Z_s&\cd     & \cd&\cd&\cd\\
  \cd&\cd&\cd   & T_t&\cd&\cd\\
  \cd&\cd&\cd  & \cd&\cd&\cd
\end{array}\right],
\end{equation}
where $X_r,\ Y_r,\ Z_s$, and\/ $T_t$
are nonsingular $r\times r$, $r\times
r$, $s\times s$, and $t\times t$ blocks
in which all diagonal entries are
nonzero and
\begin{equation}\label{3.04}
r+s=\rank A,\quad r+t=\rank B,\quad
r+s+t=\rank\: [A\,|\,B].
\end{equation}
\end{lemma}

\begin{proof}
Denote
\begin{equation*}\label{3.05}
  \rho_A=\rank A, \quad \rho_B=\rank B,\quad
\rho_M=\rank\: [A\,|\,B].
\end{equation*}

We first reduce $[A\,|\,B]$ by a
permutation of columns to the form
$[\phantom{a}\cd\ \ A_1\,|\,B\:]$,
where $A_1$ has $\rho_A$ columns and
they are linearly independent. Then we
reduce it to the form $[\phantom{a}\cd\
\ A_1\,|\,B_1\ \ \cd\phantom{a}]$,
where $[A_1\,|\,B_1]$ has $\rho_M$
columns and they are linearly
independent.

Lemma \ref{l3.1} to $[A_1\,|\,B_1]$
ensures that the matrix
$[\phantom{a}\cd\ \ A_1\,|\,B_1\ \
\cd\phantom{a}]$ reduces by a
permutation of rows to the form
\begin{equation}\label{3.06}
\begin{matrix}
  \text{$\rho_A$ rows}\: \{
  \\ {}\\ {}
\end{matrix}
\left[\begin{array}{cc|cc}
  \cd&A_2 & \cd&\cd \\
  \cd&\cd & B_2&\cd\\
    \cd&\cd & \cd&\cd
\end{array}\right]
\begin{matrix}
  \Big\}\: \text{$\rho_M$ rows}
  \\ {}
\end{matrix}
\end{equation}
with nonsingular square matrices $A_2$
and $B_2$.

Rearranging rows of the first strip and
breaking it into two substrips, we
reduce \eqref{3.06} to the form
\begin{equation}\label{3.07}
\begin{matrix}
  \text{$\rho_A$ rows}\: \Big\{
  \\ {}\\ {}
\end{matrix}
\left[\begin{array}{cc|cc}
  \cd&A_3 & \cd&\cd \\
  \cd& A_4 & B_3&\cd\\
  \cd&\cd & B_2&\cd\\
    \cd&\cd & \cd&\cd
\end{array}\right]
\begin{matrix}
 {}\\ \Big\}\: \text{$\rho_B$ rows}
  \\ {}
\end{matrix}
\begin{matrix}
 \Bigg\}\: \text{$\rho_M$ rows}
  \\ {}
\end{matrix}
\end{equation}
where the matrices
\begin{equation*}
\left[\begin{array}{cc}
  A_3 \\ \hline
  A_4
\end{array}\right]
\quad\text{ and }\quad
\left[\begin{array}{ccc}
  B_3&\cd \\ \hline
  B_2&\cd
\end{array}\right]
\end{equation*}
have linearly independent rows. Lemma
\ref{l3.1} to their transposes insures
that \eqref{3.07} reduces by
permutations of columns to the form
\begin{equation}\label{3.08}
\begin{matrix}
  \rho_A\: \Big\{
  \\ {}\\ {}
\end{matrix}
\left[\begin{array}{ccc|ccc}
  \cd&Z&\cd  & \cd&\cd&\cd \\
  \cd&\cd&X & Y&\cd&\cd\\
  \cd&\cd&\cd & \cd&T&\cd\\
    \cd&\cd&\cd & \cd&\cd&\cd
\end{array}\right]
\begin{matrix}
 {}\\ \Big\}\: \rho_B
  \\ {}
\end{matrix}
\begin{matrix}
 \Bigg\}\: \rho_M
  \\ {}
\end{matrix}
\end{equation}
with nonsingular $X,\ Y,\ Z$, and $T$.
If an $n$-by-$n$ matrix has a nonzero
determinant, then one of its $n!$
summands is nonzero, and we may dispose
the entries of this summand along the
main diagonal by a permutation of
columns. In this manner we make nonzero
the diagonal entries of $X,\ Y,\ Z$,
and $T$. At last, we reduce
\eqref{3.08} to the form \eqref{3.03}
by permutations of rows and columns.
\end{proof}

\subsection{Proof of Theorems \ref{t0.1}
and \ref{t1.2}} \label{s3}

In this section $M(x)=[A({x})\,|\,
B({x})]$ is the matrix \eqref{1.7},
${\cal A}$ and ${\cal B}$ are the
matchboxes from Theorem \ref{t1.2}, and
$r_A,\,r_B,\,r_M$ are the numbers
\eqref{1.7v}.

\begin{lemma} \label{lll}
\begin{equation}\label{e10}
\size {\cal A}= r_A,\quad \size {\cal
B}= r_B,\quad \size {{\cal A}\Cup{\cal
B}}= r_M.
\end{equation}
\end{lemma}

\begin{proof}
By Lemma \ref{l3.2}, the matrix $M(x)$
over the field $\mathbb K$ of rational
functions \eqref{3.2} reduces by
permutations of rows and by
permutations of columns within $A({x})$
and $B({x})$ to a matrix $N({x})$ of
the form \eqref{3.03}, in which by
\eqref{3.04}
\begin{equation}\label{3.3q}
r+s=r_A,\quad r+t=r_B,\quad r+s+t=r_M.
\end{equation}
The diagonal entries of $X_r,\ Y_r,\
Z_s$, and $T_t$ are all nonzero, and
hence they are independent unknowns;
replacing them by $1$ and the other
unknowns by $0$, we obtain the matrix
\begin{equation}\label{3.3aa}
N(a)=\left[\begin{array}{ccc|ccc}
  I_r&0&0     & 0&I_r&0\\
  0&I_s&0     & 0&0&0\\
  0&0&0   & I_t&0&0\\
  0&0&0  & 0&0&0
\end{array}\right],\qquad a\in
\{0,1\}^n.
\end{equation}
The inverse permutations of rows and
columns reduce $N({x})$ to $M({x})$,
and hence $N(a)$ to $M(a)$. As follows
from \eqref{3.3aa},
$a=\varepsilon_{{\cal A}'\cup{\cal
B}'}$, where ${\cal A}'$ is a left
matchbox, ${\cal B}'$ is a right
matchbox, and by \eqref{3.3q}
\[
\size {\cal A}'= r_A,\quad \size {\cal
B}'= r_B,\quad \size {{\cal
A}'\Cup{\cal B}'}= r_M.
\]
Since the matchboxes ${\cal A}$ and
${\cal B}$ are largest, $\size {\cal
A}\ge r_A$ and $\size {\cal B}\ge r_B$.
The minors $\mu_{\cal A}(x)$ of $A(x)$
and $\mu_{\cal B}(x)$ of $B(x)$
(defined in Section \ref{ss12}) are
nonzero and their orders are equal to
the sizes of ${\cal A}$ and ${\cal B}$,
hence $\size {\cal A}\le r_A$ and
$\size {\cal B}\le r_B$. We have
\[
\size {\cal A}=\size {\cal
A}'=r_A,\quad \size {\cal B}=\size
{\cal B}'=r_B,
\]
and so the matchboxes ${\cal A}'$ and
${\cal B}'$ are largest too. Because of
the minimality of the number $v({\cal
A},{\cal B})$ of common vertices and
since
\begin{equation}\label{e13}
\size {{\cal A}\Cup{\cal B}}= \size
{\cal A}+\size{\cal B}-v({\cal A},{\cal
B}),
\end{equation}
we have
\[
v({\cal A},{\cal B})\le v({\cal
A}',{\cal B}'),\quad \size {{\cal
A}\Cup{\cal B}}\ge \size {{\cal
A}'\Cup{\cal B}'}= r_M.
\]
In actuality the last inequality is an
equality since the minor $\mu_{{\cal
A}\Cup{\cal B}}(x)$ of order $r_M$ is
nonzero.
\end{proof}

\begin{lemma} \label{llll}
If $a\in{\mathbb F}^n$ and $f_{{\cal
A}{\cal B}}(a)\ne 0$, then
\begin{equation}\label{e11}
\rank A(a)=r_A,\quad \rank
B(a)=r_B,\quad \rank M(a)=r_M.
\end{equation}
\end{lemma}

\begin{proof}
The matrix $A(a)$ has a nonzero minor
$h(a)$, whose order is equal to the
rank of $A(a)$. The corresponding minor
$h(x)$ of $A(x)$ (belonging to the same
rows and columns) is a nonzero
polynomial, and so $\rank A(a)\le
\rank_{\mathbb K}A(x)=r_A$. Analogously
$\rank B(a)\le r_B$ and $\rank M(a)\le
r_M$.

By \eqref{e6}, the minors $\mu_{\cal
A}(a)$ of $A(a)$, $\mu_{\cal B}(a)$ of
$B(a)$, and $\mu_{{\cal A}\Cup{\cal
B}}(a)$ of $M(a)$ are nonzero. Their
orders are equal to the sizes of ${\cal
A},\ {\cal B}$, and ${{\cal A}\Cup{\cal
B}}$, hence
\begin{equation*}
\rank A(a)\ge \size {\cal A},\quad
\rank B(a)\ge \size {\cal B},\quad
\rank M(a)\ge \size {{\cal A}\Cup{\cal
B}}.
\end{equation*}
This proves \eqref{e11} due to
\eqref{e10}.
\end{proof}

Let $a\in{\mathbb F}^n$ and $f_{{\cal
A}{\cal B}}(a)\ne 0$. By Lemma
\ref{t1}, $M(a)$ reduces to the matrix
\eqref{1.2}, which is determined by
\eqref{1.2a} due to \eqref{1.3} and
\eqref{e11}. The matrix
$M(\varepsilon_{{\cal A}\cup{\cal B}})$
reduces by permutations of rows and
columns to the same matrix \eqref{1.2}
because \eqref{1.6x} and \eqref{e10}
imply
\begin{gather} \label{aaa}
\rank A(\varepsilon_{{\cal A}\cup{\cal
B}})=\size {\cal A}= r_A,\qquad
 \rank
B(\varepsilon_{{\cal A}\cup{\cal B}})=
\size {\cal B}= r_B,
 \\ \label{bbb}
\rank M(\varepsilon_{{\cal A}\cup{\cal
B}})= \rank M(\varepsilon_{{\cal
A}\Cup{\cal B}})= \size {{\cal
A}\Cup{\cal B}}=r_M.
\end{gather}
Hence $M(a)$ reduces to
$M(\varepsilon_{{\cal A}\Cup{\cal
B}})$. This proves Theorem \ref{t0.1}:
we can take $M_{\text{\rm gen}}$ and
$f(x)$ as indicated in \eqref{1.6dd}.
This also proves Theorem \ref{t1.2};
the equalities \eqref{1.6p} follow from
\eqref{aaa}, \eqref{bbb}, and
\eqref{e13}.
\bigskip

\noindent{\bf Acknowledgements.} Sergey
V. Savchenko read the paper and made
very important improvements and
corrections. In fact, he is a coauthor.

\end{document}